\documentclass[12pt,leqno]{article}
\newcommand{\EE}{{\rm I\kern-2pt E}}
\newcommand{\RR}{{\rm I\kern-2pt R}}
\newcommand{\DD}{{\rm I\kern-2pt D}}
\newcommand{\PP}{{\rm I\kern-2pt P}}
\newcommand{\NN}{{\rm I\kern-2pt N}}
\newcommand{\dd}{{\rm \kern 3pt I\kern-9pt d}}

\title{\large DIRICHLET FORMS IN SIMULATION}
\author{\sc Nicolas Bouleau\\
{\tt bouleau@enpc.fr}}
\date{\it ENPC, ParisTech}

\begin{document}
\maketitle

\noindent{\bf Abstract.} Equipping the probability space with a local Dirichlet form with square field operator $\Gamma$ and generator $A$
allows to improve Monte Carlo computations of expectations,  densities, and  conditional expectations, as soon as we are able to simulate a random variable $X$ 
together with $\Gamma[X]$ and $A[X]$. We give examples on the Wiener space, on the Poisson space and on the Monte Carlo space. When $X$ is real-valued
 we give an explicit formula yielding the density at the speed of the law of large numbers. \\


\noindent{1. INTRODUCTION}\\

Dirichlet forms techniques have shown their efficiency in order to obtain existence of densities under weak hypotheses, especially 
for stochastic differential equations with Lipschitz coefficients (cf [6]). We show here their utility for speeding up Monte Carlo methods
especially for the computation of such densities.

In the whole article the framework is an error structure  $(\Omega, {\cal A}, \PP,\DD,\Gamma)$, i.e. a probability space 
equipped with a local Dirichlet form with square field operator (cf [7], [4]).  We denote ${\cal E}$ the associated Dirichlet form
given by ${\cal E}[u]=\frac{1}{2}\int \Gamma[u]\,d\PP$ and $(A,{\cal D}A)$ the generator linked with ${\cal E}$ by the relation
${\cal E}[u,v]=-<A[u],v>\quad\forall u\in{\cal D}A, \forall v\in\DD.$ 

We consider a random variable $X$ belonging to the domain ${\cal D}A$ and such that we are able to simulate $X$, $\Gamma[X]$
and $A[X]$. \\

\noindent{\bf Example 1.} Wiener space.\\

Let us consider as first example a stochastic differential equation defined on the Wiener space equipped with the 
Onstein-Uhlenbeck structure (cf [12], [4], [5]) :
\begin{equation}
X_t=x_0+\int_0^t \sigma(X_s,s)dB_s+\int_0^t r(X_s,s) ds
\end{equation}
By the functional calculus for the operators
 $\Gamma$ and $A$ (cf [7], [4]), if the coefficients are smooth, the triplet $(X_t, \Gamma[X_t], A[X_t])$ is a diffusion
solution to the equation :
$$
\begin{array}{l}
\left(\begin{array}{l}
X_t\\
\Gamma[X_t]\\
A[X_t]
\end{array}\right)
=\left(\begin{array}{l}
x_0\\
0\\
0
\end{array}\right)
+{\displaystyle\int_0^t}
\left[
\begin{array}{ccc}
\sigma(X_s,s)&0&0\\
0&2\sigma^\prime_x(X_s,s)&0\\
-\frac{1}{2}\sigma(X_s,s)&\frac{1}{2}\sigma^{\prime\prime}_{x^2}(X_s,s)&\sigma'_x(X_s,s)
\end{array}
\right]
\left(\begin{array}{l}
1\\
\Gamma[X_s]\\
A[X_s]
\end{array}\right)
dB_s
\\
\\
\hspace{2.5cm}+{\displaystyle\int_0^t}
\left[
\begin{array}{ccc}
r(X_s,s)&0&0\\
\sigma^2(X_s,s)&2r^\prime_x(X_s,s)+\sigma^{\prime 2}_x(X_s,s)&0\\
0&\frac{1}{2}r^{\prime\prime}_{x^2}(X_s,s)&r^\prime_x(X_s,s)
\end{array}
\right]
\left(\begin{array}{l}
1\\
\Gamma[X_s]\\
A[X_s]
\end{array}\right)
ds
\end{array}
$$
Denoting $Y_t$ the column vector $(X_t,\Gamma[X_t],A[X_t])$ this equation writes
\begin{equation}
Y_t=Y_0+\int_0^ta(Y_s,s)dB_s+\int_0^tb(Y_s,s)ds
\end{equation}
and solving it by the Euler scheme with mesh
 $\frac{1}{n}$ on [0,T], i.e.
\begin{equation}
Y_t^n=Y_0+\int_0^ta(Y_{\frac{[ns]}{n}}^n,\frac{[ns]}{n})dB_s+\int_0^tb(Y_{\frac{[ns]}{n}}^n,\frac{[ns]}{n})ds
\end{equation}
yields a process $Y_t^n=(X_t^n, (\Gamma[X])_t^n,(A[X])_t^n)^t$. Now, it is straightforward to verify that
the second and third components of this process are  respectively equal to  $\Gamma[X_t^n]$ and $A[X_t^n]$.
In other words, if for a process  $Z$ solution to a stochastic differential equation, we denote $Z^n$ the solution to the discretized
 s.d.e. by the Euler scheme of mesh
 $\frac{1}{n}$ sur [0,T], we may write
\begin{equation}
\begin{array}{rcl}
\Gamma[X^n]&=&(\Gamma[X])^n\\
A[X^n]&=&(A[X])^n.
\end{array}
\end{equation}
Thus,  in order to compute the density of  $X_T$, if we approximate it by the Euler scheme
 $X^n_T$ and use the fact that the densities
$p_T(x_0,x)$ and $ p^n_T(x_0,x)$ of $X_T$ and $X^n_T$ are close together and satisfy under regular hypotheses
$$\sup_{x_0,x}|p_T(x_0,x)-p^n_T(x_0,x)|\leq \frac{K}{n}$$
(cf  [1], [2], [10], [11, thm 4.1]) we are eventually in a situation where we have to estimate the density of
 $X^n_T$ in a framework where we are able to simulate   $X_T^n$, $\Gamma[X_T^n]$ and $A[X_T^n]$ thanks to the relations
 (4).

\noindent{\bf Remark}. Starting from the same equation  (1), instead of putting an error on the Brownian motion, we can simply 
put an error on the initial value $x_0$. We obtain that  $(X_t,\Gamma[X_t], A[X_t])$ is still a diffusion, 
evidently different from the preceding one, but relations (4) are still valid, so that we know to simulate 
 $X_T^n$, $\Gamma[X_T^n]$ and $A[X_T^n]$.

In fact, the sequel will show that we have interest to consider both an error on the Brownian motion  and an independent 
error on the initial value because this increases 
 $\Gamma[X]$.\\

\noindent{\bf Example 2.}  Poisson space.\\

Let $(\RR^d,{\cal B}(\RR^d), \mu,d,\gamma)$ be an error structure on $\RR^d$ whose generator is denoted $(a,{\cal D}a)$ 
and let  $N$ be a Poisson point process on  $\RR^d$ with intensity measure $\mu$. The space of definition of $N$, 
$(\Omega, {\cal A}, \PP)$, may be equipped with a so-called ``white" error structure
 $(\Omega, {\cal A}, \PP,\DD,\Gamma)$  (cf [4]) with the following properties 
$$\begin{array}{ll}
\forall h\in{\cal D}a &\\
&N(h)\in\DD \;{\mbox{ and }}\;\Gamma[N(h)]=N(\gamma[h])\\
&N(h)\in{\cal D}A \;{\mbox{ and }}\;A[N(h)]=N(a[h]).
\end{array}
$$
Simulating $N(\xi)$ amounts to draw a finite (poissonian) number of independent variables with law
$\mu$, so that we are in a situation where $N(h)$, $\Gamma[N(h)]$, and $A[N(h)]$ are simulatable, the same for a regular functional 
$F(N(h_1),N(h_2),\ldots,N(h_k))$.\\

\noindent{\bf Example 3.}  Monte Carlo space.\\

Let $X$ be a random variable simulatable on the Monte Carlo space by an infinite number of calls to the random function. Let us group together
the calls with respect to which the variable 
$X$ is regular and those with respect to which it is irregular or discontinuoous (use of the rejection method, etc.) so that
$X$ may be written on the space 
 $([0,1]^{\NN},{\cal B}([0,1]^{\NN}),dx^{\NN})\times([0,1]^{\NN},{\cal B}([0,1]^{\NN}),dx^{\NN})$
$$X=F(U_0,U_1,\ldots,U_m,\ldots;V_0,V_1,\ldots,V_n,\ldots)$$
where the $U_i$'s are the coordinates of the first factor and
the $V_j$'s those of the second one, the function $F$ being regular with respect to the  $U_i$'s.

Let us put on the  $U_i$'s the error structure
$$([0,1]^{\NN},{\cal B}([0,1]^{\NN}),dx^{\NN},\DD,\Gamma)= ([0,1],{\cal B}([0,1]),dx, d,\gamma)^{\NN}
$$
where $(d,\gamma)$ is the closure of the operator $\gamma[u](x)=x^2(1-x)^2 u^{\prime2}(x)$ for $u\in{\cal C}^1([0,1]).$

If $F$ is ${\cal C}^2$ with respect to each $U_i$ sur $[0,1]$ and if the series
$$
\sum_{i=0}^\infty(\frac{1}{2}F^{\prime\prime}_{ii} U_i^2(1-U_i)^2+F^\prime_i U_i(1-U_i)(1-2U_i))
$$
converges in $L^2$, we have, $X\in{\cal D}A$ and
$$\begin{array}{rl}
A[X]=&\sum_{i=0}^\infty(\frac{1}{2}F^{\prime\prime}_{ii} U_i^2(1-U_i)^2+F^\prime_i U_i(1-U_i)(1-2U_i))\\
\Gamma[X]=& \sum_{i=0}^\infty F^{\prime2}_i U_i^2(1-U_i)^2.
\end{array}
$$
so that $X$, $\Gamma[X]$ et $A[X]$ are simulatable.\\

\noindent{2. DIMINISHING THE BIAS.}\\

Let be $(\Omega, {\cal A}, \PP,\DD,\Gamma)$  an error structure, $({\cal E},\DD)$ the associated Dirichlet form, and $(A,{\cal D}A)$
 the associated generator.

Let us explain the intuitive way. The symmetric Markov process associated to the error structure, in short time $\varepsilon$, 
induces an error on any regular random variable defined on 
 $(\Omega, {\cal A},\PP)$ whose bias is 
$\varepsilon A[X]$ and whose variance is $\varepsilon\Gamma[X]$.

Since the probability $\PP$ is invariant by the transition semi-group of the Markov process, the law of the random variable 
$X$ is nearly the same as that of
$$X+\varepsilon A[X]+\sqrt{\varepsilon}\sqrt{\Gamma[X]}\,G$$
where $G$ is an exogenous reduced Gaussian  variable  independent of ${\cal A}$. It follows first that the random variable
$X+\varepsilon A[X]$ which has the same expectation as
$X$, possesses a smaller variance than that of $X$. This is shown by the following result.

For $X\in({\cal D}A)^d$, we denote $\underline{\underline{\mbox{var}}}[X]$ the covariance matrix of 
 $X$, $A[X]$ the column vector with components 
$(A[X_1],\ldots, A[X_d])$, $\underline{\underline{\Gamma}}[X]$ the matrix $\Gamma[X_i,X_j]$ and
 $\sqrt{\underline{\underline{\Gamma}}[X]}$ the positive symmetric matrix square root of  $\underline{\underline{\Gamma}}[X]$.\\

\noindent{\bf Proposition 1.} {\it For} $X\in({\cal D}A)^d$,
$$\mbox{trace}(\underline{\underline{\mbox{var}}}[X+\varepsilon A[X]])=
\mbox{trace}(\underline{\underline{\mbox{var}}}[X])-2\varepsilon\sum_{i=1}^d{\cal E}[X_i]+\varepsilon^2\|A[X]\|^2.
$$
{\it If $A[X]$ is not zero, this quantity is minimum for $\varepsilon=\sum_i{\cal E}[X_i]/\|A[X]\|^2$ and is equal to}
$$\mbox{trace}(\underline{\underline{\mbox{var}}}[X])-2\frac{\sum_i{\cal E}[X_i]}{\|A[X]\|^2}.$$
{\bf Proof.} The result comes directly from the relation  ${\cal E}[X_i]=-<A[X_i],X_i>$.
$\hfill\diamond$\\

In order to calculate $\EE X$, it is interesting to simulate  $X+\varepsilon A[X]$ instead of
 $X$ as soon as $\varepsilon\in]0,2\sum_i{\cal E}[X_i]/\|A[X]\|^2[$.\\

We apply now the same idea to the computation of the density of 
 $X$ that we denote $f$ when it exists. Let  $g(x-m,\Sigma)$ be the density of the normal law on 
 $\RR^d$ with mean $m$ and covariance matrix $\Sigma$ supposed to be invertible.
Given $X, A[X],\underline{\underline{\Gamma}}[X]$ the conditional law of the random variable
$X+\varepsilon A[X]+\sqrt{\varepsilon}\sqrt{\underline{\underline{\Gamma}}[X]}\,G$ where $G$ 
is an independent reduced Gaussian variable, has a density 
$g(x-X-\varepsilon A[X], \varepsilon\underline{\underline{\Gamma}}[X])$. The goal  is to show, under suitable hypotheses
assuring 
$\underline{\underline{\Gamma}}[X]$ to be invertible, that $\EE g(x-X-\varepsilon A[X], \varepsilon\underline{\underline{\Gamma}}[X])$ 
converges to
$f(x)$ faster than in the classical kernel method.\\

\noindent{\bf Lemma 1. }{\it Let  $X$ be in $({\cal D}A)^d$. We suppose that $X$ possesses a conditional density 
 $\eta(x,\gamma,a)$
given $\underline{\underline{\Gamma}}[X]\!=\!\gamma$ et $A[X]\!=\!a$ such that $x\mapsto\eta(x,\gamma,a)$ be ${\cal C}^2$
with bounded derivatives. Then $\forall x\in\RR^d$}
$$
\EE[-(A[X])^t\nabla_x\eta(x,\underline{\underline{\Gamma}}[X],A[X])
+\frac{1}{2}{\mbox{trace}}\left(\underline{\underline{\Gamma}}[X].\mbox{Hess}_x\eta\right)(x,\underline{\underline{\Gamma}}[X],A[X])]=0.$$

\noindent{\bf Proof.} Let us give the argument in the case $d=1$. Let $\varphi$ be ${\cal C}^2$ with compact support
on $\RR$. By [7], denoting $A^{(1)}$ the generator in the  $L^1$ sense, we have
$$
A^{(1)}[\varphi(X)]=\varphi^\prime(X)A[X]+\frac{1}{2}\varphi^{\prime\prime}(X)\Gamma[X].
$$
Hence if $\mu(d\gamma,da)$ is the law of the pair $(\Gamma[X],A[X])$ 
$$\EE A^{(1)}[\varphi(X)]=0=\int\mu(d\gamma,da)\int(\varphi^\prime(x)a+\varphi^{\prime\prime}(x)\gamma)\eta(x,\gamma,a) dx.
$$
 Integrating by parts gives
$$\int\mu(d\gamma,da)\int\varphi(x)(-a\eta^\prime_x(x,\gamma,a)+\frac{1}{2}\gamma \eta^{\prime\prime}_{x^2}(x,\gamma,a))dx=0
$$
hence
$$\int\mu(d\gamma,da)(-a\eta^\prime_x(x,\gamma,a)+\frac{1}{2}\gamma \eta^{\prime\prime}_{x^2}(x,\gamma,a))=0
$$
as soon as $\EE|-A[X]\eta^\prime_x(X,\Gamma[X],A[X])+\frac{1}{2}\Gamma[X]\eta^{\prime\prime}_{x^2}(X,\Gamma[X],A[X])|\in L^1_{loc}(dx)
$
what is satisfied under the assumptions of the statement. \hfill$\diamond$\\

First we study the bias :\\

\noindent{\bf Proposition 2.} {\it Let  $X$ be as in the above lemma and let the conditional density
$\eta(x,\gamma, a)$
be ${\cal C}^3$ bounded with bounded derivatives.

As $\varepsilon\rightarrow 0$, the quantity
$$\frac{1}{\varepsilon^2}\left(\EE[g(x-X-\varepsilon A[X], \varepsilon\underline{\underline{\Gamma}}[X])]-f(x)\right)$$
possesses a finite limit equal to}
$$
\begin{array}{rl}
\frac{1}{2}\EE[(A[X])^t(\mbox{Hess}_x \eta)(x,\underline{\underline{\Gamma}}[X],A[X])&\!\!\!\!\!A[X]\\
-\sum_{i,j,k=1}^d&\!\!\!\!\!A[X_i]\Gamma[X_j,X_k]\eta^{\prime\prime\prime}_{x_ix_jx_k}(x,\underline{\underline{\Gamma}}[X],A[X])].
\end{array}
$$

\noindent{\bf Demonstration. } The argument begins with the relations 
$$
\begin{array}{l}
\EE[g(x-X-\varepsilon A[X], \varepsilon\underline{\underline{\Gamma}}[X])]=
\int\mu(d\gamma,da)\int g(x-y-\varepsilon, \varepsilon\gamma)
\eta(y,\gamma,a)dt \\
=\int\mu(d\gamma,da)\EE \eta(x-\varepsilon a-\sqrt{\varepsilon}\sqrt{\gamma}G,\gamma,a)
\end{array}
$$
where $G$ is a reduced Gaussian variable with values in $\RR^d$ and then consists of expanding 
$ \eta(x-\varepsilon a-\sqrt{\varepsilon}\sqrt{\gamma}G,\gamma,a)$ by the Taylor formula and taking the expectation.
The term of order zero gives $f(x)$, the term in $\sqrt{\varepsilon}$ vanishes since  $G$ is centered,
the term in $\varepsilon$ is zero because of  lemma 1, 
the term in $\varepsilon\sqrt{\varepsilon}$ vanishes because  $G^3$ is centered. The hypotheses
give the upper bounds allowing to conclude.\hfill$\diamond$\\

About the variance we have :\\

\noindent{\bf Proposition 3.} {\it Let $X$ be as in proposition 1 and such that
} $(\mbox{det}
\underline{\underline{\Gamma}}[X])^{-\frac{1}{2}}\in L^1$, 
$$
\begin{array}{l}
\lim_{\varepsilon\rightarrow 0}\varepsilon^{d/2}\EE g^2(x-X-\varepsilon A[X],\varepsilon\underline{\underline{\Gamma}}[X])
=\lim_{\varepsilon\rightarrow 0} \varepsilon^{d/2}\mbox{var}g(x-X-\varepsilon A[X],\varepsilon\underline{\underline{\Gamma}}[X])\\
=\EE\left[\frac{\eta(x,\underline{\underline{\Gamma}}[X],A[X])}{(4\pi)^{d/2}\sqrt{{\mbox{\small det}}
\underline{\underline{\Gamma}}[X]}}\right].
\end{array}
$$

\noindent{\bf Demonstration.} We have
$$\EE g^2(x-X-\varepsilon A[X],\varepsilon\underline{\underline{\Gamma}}[X])
=\int\mu(d\gamma,da)\int g^2(x-y-\varepsilon a,\varepsilon \gamma)\eta(y,\gamma,a)dy.
$$
Since
$$
g^2(z,\varepsilon\gamma)=\frac{1}{(2\pi)^{d/2}(2\varepsilon)^{d/2}\sqrt{\mbox{\small det}\gamma}}g(z,\frac{\varepsilon}{2}\gamma)
$$
we obtain the result by dominated convergence and the continuity of  $\eta$.\hfill$\diamond$\\
\newpage
\noindent{3. COMPARISON OF RATES.}\\

The quantity $\EE g(x-X-\varepsilon A[X],\varepsilon \underline{\underline{\Gamma}}[X])$ is computed by the law of large numbers
so that
the approximation $\hat{f}(x)$ of $f(x)$ is
$$
\hat{f}(x)=\frac{1}{N}\sum_{n=1}^N g(x-X_n-\varepsilon(A[X])_n,\varepsilon(\underline{\underline{\Gamma}}[X])_n)
\eqno{(*)}$$
where the indices $n$ denote independent drawings.

\noindent$\bullet$ If we are using the $L^2$ criterion
$$\|\hat{f}(x)-f(x)\|_{L^2}^2=\mbox{var}\hat{f}(x)+(\mbox{bias})^2$$
we are led to choose  $\varepsilon=\frac{1}{N^{\frac{2}{d+8}}}$ and
$$\|\hat{f}(x)-f(x)\|_{L^2}=\frac{1}{N^{\frac{4}{d+8}}}O(1)$$
to be compared with $\frac{1}{N^{\frac{2}{d+4}}}=\frac{1}{N^{\frac{4}{2d+8}}}$ in the case of the classical kernel method
 (cf [13] [14]). We see that the new method divides the dimension by two.\\

\noindent$\bullet$ The other criterion
$$c(\hat{f}(x),f(x))=\sup_{\varphi\in{\cal P}}|\EE\varphi(\hat{f}(x)-\varphi(f(x))|$$
where ${\cal P}$ is the set of polynomials of second degree  $\varphi(x)=ax^2+bx+c$ with $|a|\leq 1$ and $|b|\leq 1$, 
what gives 
$$c(\hat{f}(x),f(x))= |\EE[\hat{f}^2(x)]-f^2(x)|+|\EE\hat{f}(x)-f(x)|,$$
may be better adapted to the case of error calculus for the reason that, when the errors are thought as
germs, in short time, of Ito processes, biases have the same order of magnitude as variances (not as standard deviations).
This criterion leads us to take
$\varepsilon=\frac{1}{N^{\frac{2}{d+4}}}$ what gives $c(\hat{f}(x),f(x))=\frac{1}{N^{\frac{4}{d+4}}}O(1)$ to be compared with
 $\frac{1}{N^{\frac{2}{d+2}}}=\frac{1}{N^{\frac{4}{2d+4}}}$ in the classical case, we see that for this criterion too
the proposed method divides the dimension by two.\\

\noindent{4. DIRECT FORMULAE}\\

When  $X$ is real valued, explicit formulae may be proved that allow simulations at the speed of the law of large numbers, 
provided that in addition to $X$, $\Gamma[X]$, and $A[X]$, we are able to simulate the random variable
$\Gamma[X,\frac{1}{\Gamma[X]}]$, what is possible under additional regularity assumptions.

For instance in the  example 2, we have easily, if $X=N(h)$
$$\Gamma[X,\frac{1}{\Gamma[X]}]=\frac{N(\gamma[h,\gamma[h]])}{(N(\gamma[h]))^4}.$$

\noindent{\bf Proposition 4.} {\it a)  If $X$ is in ${\cal D}A$ with $\Gamma[X]\in\DD$ and $\Gamma[X]>0$ a.s. then $X$ has a density
$f$ possessing an l.s.c. version 
 $\tilde{f}$  and }
\begin{equation}
\tilde{f}(x)=\lim_{\varepsilon\downarrow 0}\uparrow\frac{1}{2}\EE\left(\mbox{sign}(x-X)
(\Gamma[X,\frac{1}{\varepsilon+\Gamma[X]}]+\frac{2A[X]}{\varepsilon+\Gamma[X]})\right).
\end{equation}
{\it b)
If in addition $\frac{1}{\Gamma[x]}\in\DD$,
then $X$ has a density $f$ which is absolutely continuous and }
\begin{equation}
f(x)=\frac{1}{2}\EE\left(\mbox{sign}(x-X)(\Gamma[X,\frac{1}{\Gamma[X]}]+\frac{2A[X]}{\Gamma[X]})\right).
\end{equation}

\noindent{\bf Demonstration.} Let us begin with the case {\it a)}. Since $X\in{\cal D}A$ and $\Gamma[X]\in L^2$, 
for any ${\cal C}^2$ function
$\varphi$ with bounded derivatives  (cf [7] chap I), we have $\varphi[X]\in {\cal D}A$ and
$$A[\varphi(X)]=\varphi^\prime(X)A[X]+\frac{1}{2}\varphi^{\prime\prime}(X)\Gamma[X]$$ hence $\forall\varepsilon>0$ 
\begin{equation}
\varphi^{\prime\prime}(X)=\frac{2A[\varphi(X)]+\varepsilon\varphi^{\prime\prime}(X)-2\varphi^\prime(X)A[X]}
{\varepsilon+\Gamma[X]}.
\end{equation}
Since
$\EE\frac{2A[\varphi(X)]}{\varepsilon+\Gamma[X]}=-\EE\Gamma[\varphi(X),\frac{1}{\varepsilon+\Gamma[X]}]$
taking the expectation we obtain 
\begin{equation}
\EE[\varphi^{\prime\prime}(X)\frac{\Gamma[X]}{\varepsilon+\Gamma[X]}]=
-\EE[\varphi^\prime(X)(\Gamma[X,\frac{1}{\varepsilon+\Gamma[X]}]+\frac{2A[X]}{\varepsilon+\Gamma[X]})].
\end{equation}
Let us put $K_\varepsilon(x)=\EE[\frac{\Gamma[X]}{\varepsilon+\Gamma[X]}|X=x]$
and $H_\varepsilon(x)=\EE[(\Gamma[X,\frac{1}{\varepsilon+\Gamma[X]}]+\frac{2A[X]}{\varepsilon+\Gamma[X]})|X=x]$. 
Relation
 (5) writes
\begin{equation}
\int \varphi^{\prime\prime}(x)K_\varepsilon(x)\PP_X(dx)=-\int\varphi^\prime(x)H_\varepsilon(x)\PP_X(dx).
\end{equation}
The derivative in the distributions sense of the measure  $K_\varepsilon(x)\PP_X(dx)$ is the measure
$H_\varepsilon(x)\PP_X(dx)$. It follows that the measure $K_\varepsilon(x)\PP_X(dx)$ has a density and since
$K_\varepsilon>0\;\PP_X$-a.s. (because $\Gamma[X]>0\;\PP$-a.s.) the law $\PP_X$ has a density $f$.

(We prove here again the implication $X\in\DD,\; \Gamma[X]>0\;\Rightarrow\PP_X<<dx$ which is true for any local Dirichlet form with 
square field operator cf [7])

Hence $H_\varepsilon(x)\PP_X(dx)=H_\varepsilon(x)f(x)dx$ and (9) implies that the measure$K_\varepsilon(x)\PP_X(dx)$ 
has an absolutely continuous  density
 $F_\varepsilon(x)$ and $\forall\varepsilon>0$ we have $f(x)=\frac{F_\varepsilon(x)}{K_\varepsilon(x)}$ for almost every
 $x$.

Taking $\varphi(y)=\sqrt{\lambda^2+(y-x)^2}$ in (8), it comes
$$
\EE\left(\frac{\lambda^2}{(\lambda^2+(X-x)^2)^{\frac{3}{2}}}\frac{\Gamma[X]}{\varepsilon+\Gamma[X]}\right)=
\EE\left[\frac{x-X}{\sqrt{\lambda^2+(X-x)^2}}(\Gamma[X,\frac{1}{\varepsilon+\Gamma[X]}]+\frac{2A[X]}{\varepsilon+\Gamma[X]})\right]
$$
When $\lambda\rightarrow0$, by dominated convergence, for all $x$, the right-hand side converges  to
$$
\EE[\mbox{sign}(x-X)(\Gamma[X,\frac{1}{\varepsilon+\Gamma[X]}]+\frac{2A[X]}{\varepsilon+\Gamma[X]})]
$$
where $\mbox{sign}(y)=y/|y|$ if $y=\!\!\!\!\!/\;0$ and $\mbox{sign}(0)=0$.

The left-hand side is equal to
$$\int \frac{\lambda^2}{(\lambda^2+(y-x)^2)^{\frac{3}{2}}}K_\varepsilon(y)\PP_X(dy)=
\int \frac{\lambda^2}{(\lambda^2+(y-x)^2)^{\frac{3}{2}}}F_\varepsilon(y)\,dy
$$
since $F_\varepsilon$ is continuous, this converges when $\lambda\rightarrow0$ to $2F_\varepsilon(x)$, 
therefore we have the following equality between continuous functions
$$
F_\varepsilon(x)=\frac{1}{2}\EE\left[\mbox{sign}(x-X)
(\Gamma[X,\frac{1}{\varepsilon+\Gamma[X]}]+\frac{2A[X]}{\varepsilon+\Gamma[X]})\right].
$$
Now, as $\varepsilon\downarrow0$, by its definition the function $K_\varepsilon$ increases to 1 $\;\PP_X$-a.s. since
 $\Gamma[X]$ is supposed to be strictly positive a.s. Hence $K_\varepsilon(x)f(x)$ increases to $f(x)$ 
for almost every $x$. The equality
$F_\varepsilon(x)=K_\varepsilon(x)f(x)$ valid for almost every  $x$ implies that  $F_\varepsilon$ is almost everywhere,
hence everywhere, increasing
and converges to $\tilde{f}$ l.s.c. equal to $f$ almost everywhere.

In order to prove the point  b) we proceed similarly and the hypotheses allow to replace (8) by the relation
\begin{equation}
\EE[\varphi^{\prime\prime}(X)]=
-\EE[\varphi^\prime(X)(\Gamma[X,\frac{1}{\Gamma[X]}]+\frac{2A[X]}{\Gamma[X]})].
\end{equation}
Putting $H(x)=\EE[\Gamma[X,\frac{1}{\Gamma[X]}]+\frac{2A[X]}{\Gamma[X]}|X=x]$ we see that the law of
 $X$, $\PP_X(dx)$, possesses a derivative in the sense of distributions
 $H(x)\PP_X(dx)$ which is absolutely continuous, hence$X$ has an absolutely continuous density
$f$. Taking again $\varphi(y)=\sqrt{\lambda^2+(y-x)^2}$, we obtain
$$
f(x)=\frac{1}{2}\EE\left[\mbox{sign}(x-X)(\Gamma[X,\frac{1}{\Gamma[X]}]+\frac{2A[X]}{\Gamma[X]})\right].
$$
by the same argument as above.\hfill$\diamond$\\

The density of $X$ being obtained, we can extend the formulae  (5) and  (6) in order to compute conditional expectations.\\

\noindent{\bf Proposition 5.} {\it Let be $G\in\DD\cap L^\infty$, 

\noindent a) under the assumptions of prop. 4 a), we have $dx$-a.e.}
\begin{equation}
f(x)\EE[G|X=x]=\lim_{\varepsilon\downarrow 0}\frac{1}{2}\EE\left(\mbox{sign}(x-X)
(\Gamma[X,\frac{G}{\varepsilon+\Gamma[X]}]+\frac{2GA[X]}{\varepsilon+\Gamma[X]})\right).
\end{equation}
{\it the right-hand side is l.s.c. if  $G\geq 0$,}

\noindent {\it b) under the assumptions  of prop. 4 b) and with $\frac{1}{\Gamma[x]}\in\DD\cap L^\infty$, we have 
  $dx$-a.e.}
\begin{equation}
f(x)\EE[G|X=x]=\frac{1}{2}\EE\left(\mbox{sign}(x-X)
(\Gamma[X,\frac{G}{\Gamma[X]}]+\frac{2GA[X]}{\Gamma[X]})\right).
\end{equation}
{\it where the right-hand side is continuous.}\\

The proof is similar to that of proposition 4.\hfill$\diamond$\\

Let us remark eventually that letting  $\varphi'(X)$ go to 1 in formulae (8) and (10) and in the analoguous formulae
 of proposition 5
we see that $\forall G\in\DD\cap L^\infty$
$$\EE(\Gamma[X,\frac{G}{\varepsilon+\Gamma[X]}]+\frac{2GA[X]}{\varepsilon+\Gamma[X]})=0$$ 
and also for $\varepsilon=0$
under the hypotheses of  prop. 4 b). Hence it is possible to introduce, as remarked in  [10], an arbitrary control deterministic function
$c$ in order to optimize the variance. For instance formula (12) becomes
$$
f(x)\EE[G|X=x]=\frac{1}{2}\EE\left((\mbox{sign}(x-X)-c(x))
(\Gamma[X,\frac{G}{\Gamma[X]}]+\frac{2GA[X]}{\Gamma[X]})\right).
$$

\noindent{\bf Comment.} In the kernel method  (cf [13] [12] [8]), cancelling the first term in 
the asymptotic expansion of the bias is an old idea and has been explored by several authors either by the use
of non-positive kernels (cf [9] [13]) either by a Romberg method what amounts to the preceding case.
In the method we propose in sections 2 and 3, the kernel is random and depends on the random variable
 itself. That shifts from an order of magnitude. Then the above idea may be applied again.

The nearest work to the section 4 is the study by  A. Kohatsu-Higa and R. Pettersson [10] 
which uses integration by parts on the Wiener space in the sense of Malliavin, also the paper of Bouchard, Ekeland and Touzi [3].
The difference in the  points of view comes mainly from the fact that
the integration by parts formulae are not the same, ours are simpler and do not involve Skorokhod integrals. 

Let us quote also that the results of sections 2 and 4 may be applied as well to  the Poisson space or the Monte Carlo space
with a possible choice of the Dirichlet form what gives a usefull flexibility in order to take in account the studied specific model.

\begin{list}{}
{\setlength{\itemsep}{0cm}\setlength{\leftmargin}{0.5cm}\setlength{\parsep}{0cm}\setlength{\listparindent}{-0.5cm}}
  \item\begin{center}
{\small REFERENCES}
\end{center}\vspace{0.4cm}

[1] {\sc Bally V., Talay D.} ``The law of the Euler scheme for stochastic differential equations : I. Convergence rate of the distribution function", 
{\it Prob. Th. and Rel. Fields} vol 2 No2, 93-128 (1996).

[2] {\sc Bally V., Talay D.} ``The law of the Euler scheme for stochastic differential equations : II. Convergence rate of the density", 
{\it Monte Carlo Methods and Appl.} vol 104, No1, 43-80 (1996).

[3] {\sc B. Bouchard, I. Ekeland, N. Touzi} ``On the Malliavin approach to Monte Carlo approximation of conditional expectations" {\it Finance Stochast.}
8, 45-71, (2004).

[4] {\sc Bouleau N.} {\it Error Calculus for Finance ansd Physics, the Language of Dirichlet Forms}, De Gruyter, 2003.

[5] {\sc Bouleau N.} ``Error calculus and path sensitivity in Financial models", 
{\it Mathematical Finance} vol 13/1, 115-134, Jan 2003.

[6] {\sc Bouleau N.} ``Improving Monte Carlo simulations by Dirichlet forms" C. R. Acad. Sci. Paris Ser I (2005)

[7] {\sc Bouleau N., Hirsch F.}  {\it Dirichlet Forms and Analysis on Wiener Space,} De Gruyter, 1991.

[8] {\sc M. E. Caballero, B. Fernandez, D. Nualart} ``Estimating densities and applications" {\it J. of Theoretical Probability}
 Vol 11, Nr3, (1998).

[9] {\sc Deheuvels P.} ``Estimation non param\'etrique de la densit\'e par histogrammes g\'en\'eralis\'es" 
{\it R. Statist. Appl.} Vol 25, f3, 1-24, (1977)

[10] {\sc Kohatsu-Higa A., Pettersson R.} ``Variance reduction methods for simulation of densities on Wiener space", 
{\it SIAM J. Numer. Anal.} Vol 40, No2, 431-450, (2002)

[11] {\sc Malliavin P., Thalmaier A.} ``Numerical error for SDE: Asymptotic expansion and hyperdistributions", {\it C. R. Acad. Sci. Paris}
ser. I 336 (2003) 851-856

[12] {\sc Nualart N.}  {\it The Malliavin calculus and related topics}. Springer, 1995.

[13] {\sc Parzen E.} ``On estimation of a probability density function and mode" {\it Ann. Inst. Statist. Math. }6, 127-132, (1954)

[14] {\sc Silverman B. W.} {\it Density Estimation for Statistics and Data Analysis} Chapman and Hall, 1998

\end{list}

\end{document}